\documentclass[12pt,a4paper]{article} 
\usepackage[left=1in,right=1in,top=1in,bottom=1in]{geometry} 
\usepackage{setspace} 
\usepackage{times} 
\usepackage{microtype} 

\usepackage{fancyhdr}
\usepackage{titlesec}
\titleformat{\section}{\normalfont\fontsize{12}{15}\bfseries}{\thesection}{1em}{} 
\titleformat{\subsection}{\normalfont\fontsize{12}{15}\bfseries}{\thesubsection}{1em}{} 

\usepackage{tocloft}
\usepackage{graphicx} 
\usepackage{amsmath, amssymb}
\usepackage{amsfonts}
\usepackage{amsthm}
\usepackage{graphicx}
\usepackage{hyperref}
\usepackage{algorithm, algpseudocode}
\usepackage{geometry}
\title{Generation of Random (Generalized) Orthogonal Matrices}
\author{Ali Saraeb}
\date{June 2024}

\begin{document}

\maketitle

\begin{abstract}
This paper presents an algorithmic method for generating random orthogonal matrices \(A\) that satisfy the property \(A^t S A = S\), where \(S\) is a fixed real invertible symmetric or skew-symmetric matrix. This method is significant as it generalizes the procedures for generating orthogonal matrices that fix a general fixed symmetric or skew-symmetric bilinear form. These include orthogonal matrices that fall to groups such as the symplectic group, Lorentz group, Poincaré group, and more generally the indefinite orthogonal group, to name a few. These classes of matrices play crucial roles in diverse fields such as theoretical physics, where they are used to describe symmetries and conservation laws, as well as in computational geometry, numerical analysis, and number theory, where they are integral to the study of quadratic forms and modular forms. The implementation of our algorithms can be accomplished using standard linear algebra libraries.
\end{abstract}

\tableofcontents
\newpage

\section{Introduction}

For \( N > 1 \), the matrix group \( O(N) \) is formally defined as:
\begin{align}
O(N) = \{ A \in GL_N(\mathbb{R}) : A^T A = I_N \},  
\end{align}
where \( A^T \) denotes the transpose of \( A \) and \( I_N \) represents the \( n \times n \) identity matrix. This group consists exclusively of \( n \times n \) real orthogonal matrices. It is evident that each matrix \( A \in O(N) \) possesses an inverse \( A^{-1} = A^T \). Moreover, the product of any two orthogonal matrices \( A, B \in O(N) \) remains orthogonal: \( (AB)^T = B^T A^T = B^{-1} A^{-1} = (AB)^{-1} \), confirming \( O(N) \subset GL_N(\mathbb{R}) \) as a closed subgroup under matrix multiplication, including the identity matrix \( I_N \).

The defining condition \( A^T A = I_N \) imposes \( N^2 \) scalar equations governing the \( N^2 \) real elements \( a_{ij} \) of matrix \( A \):
\begin{align}
\sum_{k=1}^{N} a_{ki} a_{kj} = \delta_{ij}, 
\end{align}
where \( \delta_{ij} \) is the Kronecker delta symbol:
\[
\delta_{ij} = \begin{cases} 
1 & \text{if } i = j, \\
0 & \text{if } i \neq j.
\end{cases}
\]

The group \( O(N) \) is known to be one of the classical compact subgroups of  \(GL_N(\mathbb{R}) \). 

The study of \( O(N) \) is motivated by their foundational role in isometries of \( \mathbb{R}^N \), preserving distances in Euclidean space. These matrix groups serve as fundamental constructs in geometry, physics, and computational sciences, providing the framework for understanding symmetry and conservation laws across diverse mathematical contexts. 

Orthogonal matrices are fundamental in the context of Euclidean spaces, arising naturally from the dot product. Specifically, an orthogonal matrix \( S \) preserves the dot product, ensuring that for any vectors \( \mathbf{u} \) and \( \mathbf{v} \) in an \( N \)-dimensional real Euclidean space, the following equality is maintained:
\[
\mathbf{u} \cdot \mathbf{v} = (S \mathbf{u}) \cdot (S \mathbf{v}).
\]

This property signifies that orthogonal matrices maintain the inner product structure of the space, thereby preserving both the lengths of vectors and the angles between them. In a broader context, consider an \( n \times n \) real invertible symmetric or skew-symmetric matrix \( S \). The corresponding \( S \)-bilinear form is defined by
\begin{align} B(\mathbf{x}, \mathbf{y}) = \mathbf{x}^{\textsf{T}} S \mathbf{y} = \sum_{i,j=1}^{N} x_i S_{ij} y_j \end{align} 
for vectors \( \mathbf{x}, \mathbf{y} \in \mathbb{R}^N \). Similarly, it is of interest to consider the linear transformations, $A,$ that preserve the S-bilinear form; that is, \begin{align}
    \mathbf{x}^{\textsf{T}} S \mathbf{y}= B(\mathbf{x}, \mathbf{y})= B(A\mathbf{x}, A \mathbf{y})= \mathbf{x}^{\textsf{T}} ( A^{\textsf{T}} S A ) \mathbf{y} \text{ for all } \mathbf{x}, \mathbf{y} \in \mathbb{R}^N.\label{(4)} \end{align} 
    
    This introduces an analogous structure known as the generalized orthogonal group of S-orthogonal matrices:
\begin{align} O_S(N) = \{ A \in GL_N(\mathbb{R}) : A^{\textsf{T}} S A = S \}. \end{align} 

This group is a closed subgroup of \( GL_N(\mathbb{R}) \), and since \( \det(S) \neq 0 \), every matrix \( A \in O(S) \) satisfies \( \det(A) = \pm 1 \). In analogy with the case of orthogonal matrices (corresponding to \( S = I_N \)) where generation of isometries is desirable (see \cite{Mezzadri} and \cite{Genz1998}), a natural question arises: how can we generate linear transformations that preserve the \( S \)-bilinear form, i.e., linear transformations that satisfy equation \ref{(4)}?
This paper addresses this question by presenting an algorithm for generating random orthogonal matrices from $O_S(N) \cap O(N)$. We now illustrate two fundamental classes of matrices that serve as simple applications of our algorithms.
\subsection{Symplectic Matrices} 

Symplectic matrices, pivotal in both mathematics and quantum physics, are \( 2n \times 2n \) matrices satisfying \( A^T \Omega A = \Omega \), where \( \Omega = \begin{bmatrix} 0 & I_N \\ -I_N & 0 \end{bmatrix} \). This property ensures preservation of the symplectic form, essential in classical mechanics for conserving phase space area in Hamiltonian systems. In quantum physics, symplectic matrices are fundamental in describing Gaussian transformations of quantum states, crucial for manipulating quantum information. By setting \( S = \Omega \), our algorithm allows for the generation of orthogonal symplectic matrices.
\subsection{Indefinite Orthogonal Matrices}
The indefinite orthogonal group \( O(p, q) \) constitutes a fundamental Lie group of linear transformations, on an \( N \)-dimensional real vector spaces, that preserve a non-degenerate symmetric bilinear form of signature \((p, q)\), where \( N = p + q \). The group's structure is uniquely determined by the signature of its bilinear form. In the realm of real vector spaces, \( O(p, q) \) is defined to be the group of matrices that preserve a specific bilinear form defined by the diagonal matrix \( g \):
\[
g = \text{diag}( \underbrace{1, \ldots, 1}_p, \underbrace{-1, \ldots, -1}_q).
\]

This bilinear form, \( [x, y]_{p, S} = \langle x, gy \rangle \), is preserved by a matrix \(A \in O(p, q) \)  if and only if \( g A^{\text{T}} g = A^{-1} \). For example, the Lorentz group, denoted as \( \text{Lor} = O(1, 3) \), represents the fundamental symmetry of spacetime in physics. It describes transformations that preserve the quadratic form \( t^2 - x^2 - y^2 - z^2 \) on \( \mathbb{R}^4 \), distinguishing time from spatial dimensions with a signature of \( (1, 3) \).  By setting \( S = g \), our algorithm allows for the generation of orthogonal indefinite orthogonal matrices.

\section{Numerical Linear Algebra Tools}

In practice, generating elements of \( O_S \) is particularly interesting when the matrix \( S \) is real invertible symmetric or skew-symmetric. Thus, we present our algorithm for these cases. The tools utilized in this paper are primarily from advanced numerical linear algebra. The case where \( S \) is symmetric is easier to handle, while the other case requires specific treatments, as demonstrated below in section $3$. We outline the main tools employed in this paper and provide sketches of their proofs.

\textbf{Theorem 1.}\textbf{ [The QR Factorization $\cite{NumericalLinearAlgebra} \relax$]}

If \( A \) is an \( N \times N \) complex square matrix, then there is a decomposition \( A = QR \) where \( Q \) is a unitary matrix and \( R \) is an upper triangular matrix. If A is a real matrix, \( Q \) will be an orthogonal matrix.

\textbf{Remark 1.1.} If \( A \) is invertible, then the factorization is unique if we require the diagonal elements of \( R \) to be positive. \\
\textit{Sketch of the proof of Theorem 1.}
Apply the Gram-Schmidt process to the columns of \( A = [ a_1| \ldots| a_N ] \). Suppose \( k-1 \) orthonormal vectors \( \{q_1, \ldots, q_{k-1}\} \) have been constructed from the \( k-1 \) columns of \( A \) (denote their span as \( S_{k-1} \)). Let the \( k \)-th column of \( Q \), $q_k$, be the orthonormal vector given by the Gram-Schmidt process:
\begin{align*}
u_k = a_k - \sum_{i=1}^{k-1} \text{proj}_{q_i} a_k = a_k - \sum_{i=1}^{k-1} \langle a_k, q_i \rangle q_i,
\end{align*}
and 
\begin{align*}
q_k = \frac{u_k}{\|u_k\|_2} \quad \text{if} \quad u_k \neq 0.
\end{align*}
If \( u_k = 0 \), make a choice of \( q_k \perp S_{k-1} \). Next, for fixed \( k \) and \( 1 \leq i \leq k-1 \), let \( R = (r_{ik}) \) where \( r_{ik} = \langle a_k, q_i \rangle \) and the remaining entries of the $k$-th column of $R$ are set to zero. This completes the proof. \(\blacksquare\) 

\textbf{Remark 1.2.} If \( A \) is invertible, the vectors \( u_k \) from the proof will never be zero.

\textbf{Theorem 2. [The Complex Schur Decomposition $\cite{MatrixComputations}] \relax$ }

If \( M \) is an \( N \times N \) square matrix with complex entries, then \( M \) can be expressed as:
\[
M = UTU^{*}
\]
for some unitary matrix \( U \) and some upper triangular matrix \( T \). This is known as the Schur form of \( M \).

\textbf{Remark 2.1.} If \( M \) is symmetric, then \( T \) is a real diagonal matrix of the eigenvalues of \( M \), and \( U \) is real.

\textbf{Remark 2.2.} The Schur decomposition extends the eigenvalue decomposition for symmetric matrices, differing by selecting an orthonormal basis of eigenvectors.

\textit{Sketch of the proof of Theorem 2.}
The characteristic polynomial of \( M \) is of degree \( N \) and has \( N \) roots by the fundamental theorem of algebra. Choose an eigenvalue \(\lambda\) (a root of the polynomial) and a corresponding normalized eigenvector \( q_1 \), so that \( M q_1 = \lambda q_1 \) and \( \|q_1\|_2 = 1 \). Let \( Q_1 = [q_1 \, | \, Q_2] \), where \( Q_2 \) is an \( n \times (n-1) \) matrix with columns chosen such that \( q_1 \) and the columns of \( Q_2 \) form an orthonormal basis of \( \mathbb{C}^n \) (achievable via the Gram-Schmidt process). Then,
\[
Q_1^* M Q_1 = 
\begin{bmatrix}
q_1^* M q_1 & q_1^* M Q_2 \\
Q_2^* M q_1 & Q_2^* M Q_2
\end{bmatrix} = 
\begin{bmatrix}
\lambda q_1^* q_1 & q_1^* M Q_2 \\
\lambda Q_2^* q_1 & Q_2^* M Q_2
\end{bmatrix} = 
\begin{bmatrix}
\lambda & q_1^* M Q_2 \\
0 & N
\end{bmatrix},
\]
where the latter is evident since \( q_1 \) is normalized and orthogonal to the columns of \( Q_2 \), and \( N = Q_2^* M Q_2 \). Repeat the process on \( N \) to complete the proof. \(\blacksquare\)

\textbf{Theorem 3. [The Real Schur Decomposition $\cite{MatrixComputations} \relax$] }

If \( A \) is a real \( N \times N \) matrix then $A$ has a real Schur decomposition \( A = Q T Q^T \) where \( Q \) is an orthogonal matrix and \( T \) is a real matrix of the form

\[
T = \begin{bmatrix}
T_1 & \times & \cdots & \times \\
O & T_2 & \times & \times \\
O & O & \ddots & \times \\
O & O & O & T_j
\end{bmatrix}
\]
\noindent 
where the \( T_j \)'s are either \( 1 \times 1 \) or \( 2 \times 2 \) matrices. Each \( 2 \times 2 \) block will correspond to a pair of complex conjugate eigenvalues of \( A \). Moreover, if \( A \) is also symmetric, then the above decomposition coincides with the spectral decomposition with the additional property that \( Q \) is orthogonal. Also, if \( A \) is an invertible even-sized skew-symmetric matrix with its pure imaginary pairs of eigenvalues \(\{ (-\lambda_k i,\lambda_k i)\}_{k=1}^{n/2}\) (\(\lambda_k > 0\) without loss of generality), then \( T \) is an invertible quasi-diagonal skew-symmetric matrix, and the \( T_j \)'s are \( 2 \times 2 \) matrices of the form 
\begin{align*}
T_j =
\begin{bmatrix}
    0 & \lambda_{k_j} \\
    -\lambda_{k_j} &  0
\end{bmatrix} \quad \text{for some} ~ 1 \le k_j \le n/2.
\end{align*}

\textbf{Remark 3.1.} There is no odd-sized invertible skew-symmetric matrix.

\textbf{Remark 3.2.} If \( A \) is an odd-sized skew-symmetric matrix, then at least one of the \( T_j \)'s is a \( 1 \times 1 \) zero matrix.

\textit{Proof of Theorem 3.}
One proceeds by induction just as in the case of the complex Schur decomposition. The proof depends on the property that for each pair of complex conjugate eigenvalues of \( A \) there is a 2-dimensional subspace of \( \mathbb{R}^n \) that is invariant under \( A \) (see $ \cite{MatrixComputations}\relax$).  \(\blacksquare\)

\section{Generation of Matrices from  \texorpdfstring{\(O_S(N) \cap O(N)\)}{TEXT}}

In this section, we present and detail our method for generating matrices from the group \(O_S(N) \cap O(N)\), under two scenarios: \(S\) being a real invertible symmetric matrix and \(S\) being a real invertible skew-symmetric matrix. It is noteworthy that the case when \(S\) is symmetric is comparatively simpler than the case when \(S\) is skew-symmetric. We shall first address the simpler scenario.
\subsection{Case 1: \texorpdfstring{$S$}{TEXT}  is Symmetric}

If \(S\) is an $N \times N$ real symmetric matrix, it has only real eigenvalues. In this case, we can fix a choice of the real Schur decomposition of $S$ and write
\begin{equation}
    S = U T U^T,
\end{equation}
where \(U\) is an orthogonal matrix and \(T\) is a real diagonal matrix. 

The problem of generating an element of \(O_S(N) \cap O(N)\) reduces to generating an orthogonal matrix \(A\) satisfying 
\begin{equation}
    A^T U T U^T A = U T U^T,
\end{equation}
which is equivalent to 
\begin{equation}
    (U^T A^T U) T (U^T A U) = T.
\end{equation}
Thus, it suffices to generate an orthogonal matrix \(B\) satisfying 
\begin{equation}
    B^T T B = T,
\end{equation}
since we can then set \(A = U B U^T\). \\ Now we can write  \begin{align}
T = \begin{bmatrix}
\lambda_1 I_{k_1} & 0 & \cdots & 0 \\
0 & \lambda_2 I_{k_2} & \cdots & 0 \\
\vdots & \vdots & \ddots & \vdots \\
0 & 0 & \cdots & \lambda_m I_{k_m}
\end{bmatrix},
\end{align}
where \(\{\lambda_i\}_i\) are the distinct eigenvalues of $T$, $k_i$ is the multiplicity of $\lambda_i$, and \(I_{k_i}\) is the identity matrix of size \(k_i \times k_i\). We now present an elementary fact from linear algebra. 

\textbf{Lemma.} An orthogonal matrix \(B\) commutes with a diagonal matrix \(T\) of the form $(10)$
if and only if \(B\) has the form
\[
B = \begin{bmatrix}
B_1 & 0 & \cdots & 0 \\
0 & B_2 & 0 & 0 \\
0 & 0 & \ddots & 0 \\
0 & 0 & 0 & B_m
\end{bmatrix},
\]
where \(B_i\) is a \(k_i \times k_i\) orthogonal matrix for all \(i\).

\textit{Proof of the Lemma.} The backward direction is immediate. Suppose \(B\) is an orthogonal matrix that commutes with \(T\), i.e., \(BT = TB\).
Assume \(B\) is of the form
\[
B = \begin{bmatrix}
B_{11} & B_{12} & \cdots & B_{1m} \\
B_{21} & B_{22} & \cdots & B_{2m} \\
\vdots & \vdots & \ddots & \vdots \\
B_{m1} & B_{m2} & \cdots & B_{mm}
\end{bmatrix},
\]
where each \(B_{ij}\) is a \(k_i \times k_j\) block.

The commutation condition \(BT = TB\) gives
\[
\begin{bmatrix}
\lambda_1 B_{11} & \lambda_2 B_{12} & \cdots & \lambda_m B_{1m} \\
\lambda_1 B_{21} & \lambda_2 B_{22} & \cdots & \lambda_m B_{2m} \\
\vdots & \vdots & \ddots & \vdots \\
\lambda_1 B_{m1} & \lambda_2 B_{m2} & \cdots & \lambda_m B_{mm}
\end{bmatrix}
=
\begin{bmatrix}
\lambda_1 B_{11} & \lambda_1 B_{12} & \cdots & \lambda_1 B_{1m} \\
\lambda_2 B_{21} & \lambda_2 B_{22} & \cdots & \lambda_2 B_{2m} \\
\vdots & \vdots & \ddots & \vdots \\
\lambda_m B_{m1} & \lambda_m B_{m2} & \cdots & \lambda_m B_{mm}
\end{bmatrix}.
\]

For the latter equality to hold, each off-diagonal block \(B_{ij}\) must be zero because \(\lambda_i \neq \lambda_j\) when \(i \neq j\). Therefore, \(B\) must be block diagonal, and since \(B\) is orthogonal, each \(B_{ii}\) must be orthogonal. \(\blacksquare\) \\

Thus, generating an element $A$ of \(O_S(N) \cap O(N)\) is equivalent to generating a block diagonal orthogonal matrix \(B\). Since for all \(i \neq j\), every entry of a block matrix \(B_i\) is pairwise independent and identically distributed with the entries of a block matrix \(B_j\) due to the block matrix multiplication and the orthogonality definition $B^{T} B= I_n$, one can generate $A$ by the following steps:

\begin{enumerate}
    \item For each \(i\), generate a random matrix \(A_i \in \mathrm{GL}(k_i, \mathbb{R})\), where the entries of \(A_i\) are independently and identically distributed standard real random variables.
    \item For each $i$, extract a random \(k_i \times k_i\) orthogonal matrix \(B_i\) by applying the QR decomposition to \(A_i\). To ensure that the matrices \(B_i\) (and hence \(B\)) are suitably random, apply the QR decomposition in such a way that it is unique. Specifically, enforce the diagonal entries of the upper triangular matrices \(R_i\) obtained from the QR decomposition of \(A_i\) to be positive (see \cite{Mezzadri}). An alternative approach to generating the random orthogonal matrices \( B_i \) is to employ the polar decomposition. When applied to an invertible matrix, the polar decomposition is unique, thereby ensuring that the orthogonal matrices are appropriately generated.
    \item Form the block diagonal matrix \(B = \text{diag}(B_1, B_2, \ldots, B_m)\).
    \item Finally, return \(A = U B U^T\)
\end{enumerate}

\begin{algorithm}
\caption{Generation of matrices from \(O_S(N) \cap O(N)\), where $S$ is invertible symmetric}
\begin{algorithmic}[1]
\State \textbf{Input:} An invertible symmetric $N \times N$ matrix \(S\)
\State \textbf{Output:} An orthogonal $S$-orthogonal matrix \( A \)
\For{$i = 1$ to $m$}
    \State $A_i =rand(k_i)$.
    \State Compute the QR decomposition of $A_i$: $A_i = Q_i R_i$
    \State Ensure uniqueness:
        \State \quad Compute $\Lambda = \text{diag}\left(\frac{r_{ii}}{|r_{ii}|}\right)$, where $r_{ii}$ are the diagonal entries of $R_i$.
        \State \quad Set $Q_i' = Q_i \Lambda$, ensuring $\Lambda^{-1} R_i$ has positive diagonal entries.
    \State Set $B_i = Q_i'$
\EndFor
\State Form the block diagonal matrix $B = \text{diag}(B_1, B_2, \ldots, B_m)$
\State Compute $A = U B U^T$
\State \Return $A$
\end{algorithmic}
\end{algorithm}

\subsection{Case 2: \texorpdfstring{$S$}{TEXT} is Skew-Symmetric}

For a \( 2N \times 2N \) skew-symmetric invertible matrix \( S \), the spectrum \( \sigma \) of $S$ consists of pairs of complex conjugates of purely imaginary complex numbers. This scenario requires a more intricate approach since the spectral decomposition or the complex Schur decomposition factors \( S \) into complex matrices. Thus, in this case, we apply the real Schur decomposition \( S = U T U^T \), where \( U \) is an orthogonal matrix and \( T \) is quasi-diagonal skew-symmetric invertible matrix of the form
\begin{align}
T &= \begin{bmatrix}
T_1 & 0 & \cdots & 0 \\
0 & T_2 & 0 & 0 \\
0 & 0 & \ddots & 0 \\
0 & 0 & 0 & T_m
\end{bmatrix}, \\
\intertext{where the \( T_j \)'s are \( 2 \times 2 \) matrices of the form}
T_j &= \begin{bmatrix}
    0 & \lambda_{k_j} \\
    -\lambda_{k_j} &  0
\end{bmatrix} \quad \text{for some } 1 \le k_j \le N.
\end{align}
The numbers \( \lambda_{k_j} \) are the eigenvalues of \( T \) with positive imaginary parts, ~as described ~in Theorem 3.

Given this setup, the problem reduces to generating an orthogonal matrix \( A \) such that 
\begin{align}
A^T U T U^T A = U T U^T,
\end{align}
which simplifies to 
\begin{align}
(U^T A^T U) T (U^T A U) = T.
\end{align}
Therefore, it suffices to generate an orthogonal matrix \( B \) satisfying 
\begin{align}
B^T T B = T,
\end{align}
and then set \( A = U B U^T \).

Finally, we may permute the rows and columns of \( T \) by an orthogonal permutation matrix \( P \) so that \( J := P T P^T = \begin{bmatrix}
    0_N & D \\
    -D & 0_N
\end{bmatrix} \), where \( D = \text{diag}(\lambda_1, \lambda_2, \ldots, \lambda_N) \). Hence, generating a matrix \( B \) that satisfies equation (13) is now equivalent to generating a matrix \( C \) such that 
\begin{align}
C^T J C = J,
\end{align}
as then we can set \( B = P^T C P \).

We next provide a characterization of the class of orthogonal matrices satisfying equation (16). In other words, we characterize the groups \( O_J(2N) \cap O(2N) \), where \( J \) has the above form. The characterization is summarized in the following proposition that generalizes a well-known fact about symplectic matrices (see [\cite{deGosson}, p. 16]).

\textbf{Definition} For an invertible diagonal matrix \( D \), let \( U_D(N, \mathbb{C}) \) be the subgroup of the group \( U(N, \mathbb{C}) \) of \( N \times N \) unitary matrices that commute with the diagonal matrix \( D \). In other words, \( U_D(N, \mathbb{C}) := \{ Q \in U(N, \mathbb{C}) : Q D = D Q \} \). Note that \( U(N, \mathbb{C}) = U_{I_N}(N, \mathbb{C}) \), where \( I_N \) is the \( N \times N \) identity matrix.

\textbf{Proposition}
Let \( J = \begin{bmatrix}
    0_N & D \\
    -D & 0_N
\end{bmatrix} \), where \( D \) is an invertible diagonal matrix. The monomorphism \( \mu : M(N, \mathbb{C}) \to M(2N, \mathbb{R}) \) defined by
\begin{align}
u = A + iB \mapsto \mu(u) = \begin{bmatrix} A & -B \\ B & A \end{bmatrix}
\end{align}
(for \( A \) and \( B \) real) identifies the group \( U_D(N, \mathbb{C}) \) with the subgroup
\begin{align}
U(N) = O_J(2N) \cap O(2N)
\end{align}
of \( O_J(2N) \).

\textit{Proof of the Proposition.}
It suffices to prove the equality in equation (18). We begin by making the following observations. First, an \( N \times N \) complex matrix \( u = A + iB \) (\( A, B \) real matrices) is unitary iff \( u^{*} u = I_N \) iff \( A^t A + B^t B + i(A^t B - B^t A) = I_N \) iff \( A^t A + B^t B = I_N \) and \( A^t B = B^t A \). The latter is equivalent to the conditions \( A A^t + B B^t = I_N \) and \( B A^t = A B^t \). Second, an \( N \times N \) complex matrix \( u = A + iB \) (\( A, B \) real matrices) commutes with \( D \) iff \( A \) and \( B \) commute with \( D \). Third, if \( Q = \begin{bmatrix} A & -B \\ B & A \end{bmatrix} \), where \( A \) and \( B \) are real \( N \times N \) matrices and \( u = A + iB \) commutes with \( D \), then
\begin{align}
Q^t J Q = \begin{bmatrix} A^t D B - B^t D A & A^t D A + B^t D B \\ -(A^t D A + B^t D B) & A^t D B - B^t D A \end{bmatrix} = \begin{bmatrix} 0_N & D \\ -D & 0_N \end{bmatrix},
\end{align}
where the last equality follows from the conditions in the first and second observations.

We now return to prove the equality in equation (18). If \( Q = \mu(u) \), where \( u \in U_D(N, \mathbb{C}) < U(N, \mathbb{C}) \), then it follows by block-wise matrix multiplication and the first observation that
\begin{align}
Q^{-1} = \begin{bmatrix} A^T & B^T \\ -B^T & A^T \end{bmatrix} = Q^T,
\end{align}
hence \( Q \in O(2n) \). Moreover, it follows by the third observation that \( Q \in O_J(2N) \), which proves the inclusion \( U(n) \subset Sp(n) \cap O(2n) \).

Conversely, suppose that \( Q \in O_J(2N) \cap O(2n) \). Then
\begin{align}
JQ = (Q^T)^{-1} J = QJ,
\end{align}
which implies that \( Q \in U(n) \), so \( O_J(2N) \cap O(2n) \subset U(n) \). \( \blacksquare \)

Thus, since it suffices to generate a \(2N \times 2N\) orthogonal matrix \(C\) satisfying equation (16), it suffices to generate a unitary matrix \(Q\) that commutes with the diagonal matrix \(D = \text{diag}(\lambda_1, \ldots, \lambda_N)\). In other words, we would like to generate an \(N \times N\) unitary matrix \(Q\) such that
\begin{align}
Q^{*} D Q = D.
\end{align}
Again, we may introduce a unitary permutation matrix \(W\) such that 
\begin{align}
 D' := W D W^{*} = W D W^{t} = \begin{bmatrix}
\mu_{1} I_{k_1} & 0 & \cdots & 0 \\
0 & \mu_{2} I_{k_2} & \cdots & 0 \\
\vdots & \vdots & \ddots & \vdots \\
0 & 0 & \cdots & \mu_{M} I_{k_M}
\end{bmatrix},
\end{align} where \(\mu_{j}'s\) are the distinct eigenvalues of \(D\) with respective multiplicities \(k_j\). Thus, it suffices to generate an \(N \times N\) unitary matrix \(V\) such that
\begin{align}
V^{*} D' V = D',
\end{align}
as then we can set \(Q = W^{*} V W\). However, just as in case 1, we have the following elementary linear algebra fact.

\textbf{Lemma.} A unitary matrix \(V\) commutes with a diagonal matrix \(D'\) of the form
\begin{align}
D' = \begin{bmatrix}
\mu_{1} I_{k_1} & 0 & \cdots & 0 \\
0 & \mu_{2} I_{k_2} & \cdots & 0 \\
\vdots & \vdots & \ddots & \vdots \\
0 & 0 & \cdots & \mu_{M} I_{k_M}
\end{bmatrix}
\end{align}
if and only if \(V\) has the block diagonal form
\begin{align}
V = \begin{bmatrix}
V_1 & 0 & \cdots & 0 \\
0 & V_2 & 0 & 0 \\
0 & 0 & \ddots & 0 \\
0 & 0 & 0 & V_m
\end{bmatrix},
\end{align}
where \(V_j\) is a \(k_j \times k_j\) unitary matrix for all \(j\). \\

Thus, similar to case 1 in section 3.2, since for all \(i \neq j\), every entry of a block matrix \(B_i\) is pairwise independent and identically distributed with the entries of a block matrix \(B_j\), one can apply the QR factorization to random invertible complex matrices as described in Algorithm $1$ to generate the desired matrix (see Algorithm $2$ on page $12$).
\begin{algorithm}
\caption{Generation of matrices from \(O_S(2N) \cap O(2N)\), where $S$ is invertible skew-symmetric}
\begin{algorithmic}[1]
\State \textbf{Input:} A \( 2N \times 2N \) invertible skew-symmetric matrix \( S \)
\State \textbf{Output:} An orthogonal \( S \)-orthogonal matrix \( A \)
\State Perform real Schur decomposition: \( S = U T U^T \)
\State Use a permutation matrix \( P \) to obtain \( J := P T P^T = \begin{bmatrix} 0_N & D \\ -D & 0_N \end{bmatrix} \), where \( D \) is a real diagonal matrix.
\State Use a permutation matrix \( W \) to obtain \( D' := W D W^* = \operatorname{diag}( \mu_{1} I_{k_1}, \mu_{2} I_{k_2}, \ldots, \mu_{M} I_{k_M}) \), where $\{\mu_{j}\}_{j=1}^M$ are the distinct eigenvalues of \( D \).
\State Generate a unitary matrix \( Q \) such that \( Q^* D Q = D \)
    \State \quad Generate a unitary matrix \( V \) such that \( V^* D' V = D' \)
    \State \quad \quad \textbf{for} $i = 1$ to $M$ \textbf{do}
        \State \quad \quad \quad Generate a random complex matrix \( V_i = \operatorname{rand}(k_i) + \operatorname{rand}(k_i) \cdot i \)
        \State \quad \quad \quad Compute the QR decomposition of \( V_i \): \( V_i = Q_i R_i \)
        \State \quad \quad \quad \quad Ensure uniqueness:
            \State \quad \quad \quad  \quad \quad Compute \( \Lambda = \operatorname{diag}\left(\frac{r_{ii}}{|r_{ii}|}\right) \), where \( r_{ii} \) are the diagonal entries of \( R_i \)
            \State \quad \quad \quad \quad \quad Set \( Q_i' = Q_i \Lambda \), ensuring \( \Lambda^{-1} R_i \) has positive diagonal entries
        \State \quad \quad \quad Set \( V_i = Q_i' \)
    \State \quad \quad \textbf{end for}
    \State \quad \quad  Form the block diagonal matrix \( V = \operatorname{diag}(V_1, V_2, \ldots, V_M) \)
    \State \quad Set \( Q = W^* V W \)
\State Generate a \( 2N \times 2N \) orthogonal matrix \( C \) such that \( C^T J C = J \) by setting \( C = \mu(Q) \)
\State Compute \( B = P^T C P \)
\State Compute \( A = U B U^T \)
\State \Return \( A \)
\end{algorithmic}
\end{algorithm}
\newpage

\section*{Acknowledgment} The author would like to express gratitude to Professors Ghaith Hiary and Richard Aoun for their assistance in answering specific queries during the course of this research.


\begin{thebibliography}{99} 
\bibitem{MatrixComputations}
Gene H. Golub and Charles F. Van Loan,
\textit{Matrix Computations},
3rd edition,
Johns Hopkins University Press,
1996,
ISBN: 978-0801854149.
\bibitem{NumericalLinearAlgebra}
Lloyd N. Trefethen and David Bau III,
\textit{Numerical Linear Algebra},
SIAM,
1997,
ISBN: 978-0898713619.

\bibitem{Mezzadri}
Francesco Mezzadri,
\textit{How to Generate Random Matrices from the Classical Compact Groups},
Notices of the American Mathematical Society,
vol. 54, no. 5,
2006,
pp. 592–604.

\bibitem{Genz1998}
    Alan Genz,
    \texorpdfstring{\textit{Methods for generating random orthogonal matrices}, Monte-Carlo and Quasi-Monte Carlo Methods 1998: Proceedings of a Conference held at the Claremont Graduate University, Claremont, California, USA, June 22--26, 1998, Springer Berlin Heidelberg, 2000.}{Methods for generating random orthogonal matrices}

\bibitem{deGosson}
Maurice de Gosson,
\textit{Introduction to Symplectic Mechanics: Lectures I-II-III},
World Scientific, Singapore, 2015.


\end{thebibliography}
\end{document}